\documentclass[12pt]{amsart}

\usepackage{graphicx,amssymb,xcolor,enumerate,latexsym}
\usepackage{hyperref}
\hypersetup{
bookmarks=true,         %show bookmarks bar?
pdffitwindow=false,     % window fit to page when opened
pdfstartview={FitH},    % fits the width of the page to the window
colorlinks=true,      %false: boxed links; true: colored links
citecolor=red,
}

\usepackage{enumitem}

\usepackage{float}

\usepackage{mathrsfs} 
\usepackage[top=3cm, bottom=3cm, left=2.4cm, right=2.4cm]{geometry} 

\usepackage{palatino} %palatino for English
\usepackage{mathpazo} %palatino for math formulas
\usepackage{amsmath} % ams for math formula
\usepackage{bm}

\makeatletter
\def\l@subsection{\@tocline{2}{0pt}{2.5pc}{5pc}{}}
\makeatother

\DeclareSymbolFont{largesymbol}{OMX}{yhex}{m}{n}
\DeclareMathAccent{\Widehat}{\mathord}{largesymbol}{"62}
\newcommand*\di{\mathop{}\!\mathrm{d}}

\numberwithin{equation}{section}              % numerazionedelleequazioni
\newtheorem{theorem}{Theorem}[section]

\newtheorem{lemma}{Lemma}[section]
\newtheorem{proposition}{Proposition}[section]
\newtheorem*{proposition*}{Proposition}

\newtheorem*{corollary*}{Corollary}
\newtheorem{definition}{Definition}[section]
\newtheorem*{definitions*}{Definitions}

\newtheorem*{acknowledgements*}{Acknowledgements}

\newtheorem*{conjecture*}{\bf Conjecture}

\newtheorem*{example*}{\bf Example}
\theoremstyle{remark}
\newtheorem{remark}{\bf Remark}[section]

\usepackage{amsthm}
\usepackage{amssymb}

\newenvironment{proof of mainresult1}[1][Proof]{\proof[\textbf{Proof of Thoerem \ref{thm:globalexistence}}]}{\endproof}

\newenvironment{proof of mainresult2}[1][Proof]{\proof[\textbf{Proof of Thoerem \ref{thm:convergence}}]}{\endproof}

\usepackage{setspace}

\begin{document}

\author{Cong Wang}
\address[C. Wang]{School of Mathematical Sciences, Fudan University,  Shanghai,	200433, P. R. China.}
\email{math\_congwang@163.com}

\title[Fokker-Planck equation on locally finite graphs]{Gradient flow structure, well-posedness and asymptotic behavior of Fokker-Planck equation on locally finite graphs}

\begin{abstract}
This paper investigates the gradient flow structure, well-posedness, and asymptotic behavior of the Fokker-Planck equation defined on locally uniformly finite graphs, which is highly non-trivial compared with the finite case. We first construct a 2-Wasserstein-type metric and gradient flow equation in the probability density space associated with the underlying graphs. Then, we prove the global existence of solution to the Fokker-Planck equation using a novel approach that differs significantly from the methods applied in the finite case. We also demonstrate that the solution converges to the Gibbs distribution in the $\ell^{r}(V,\bm{\pi})$ norm with $r\in [2,\infty]$, by using the indicator set partitioning method. To the best of our knowledge, this work seems the first result on the study of Wasserstein-type metrics and the Fokker-Planck equation in probability density spaces defined on infinite graphs.
\end{abstract}

\maketitle

\section{Introduction}\label{sec:into}
The classical Fokker-Planck equation describes the evolution of the probability density for a stochastic process associated with an It\^{o} stochastic differential equation. The seminal work by Jordan, Kinderlehrer, and Otto \cite{jordan1998sima} revealed the connection among the Wasserstein metric (also named the Monge-Kantorovich metric), the Fokker-Planck equation, and the associated free energy functional, which is a linear combination of a potential energy functional and the negative of the Gibbs-Boltzmann entropy functional. In fact, the Fokker-Planck equation can be interpreted as a gradient flow, or a steepest descent, for the free energy with respect to the 2-Wasserstein metric. This discovery has served as a starting point for numerous developments in evolution equations, probability theory, and geometry \cite{ambrosio2006springer,villani2003ams,villani2008springer}. 

In recent years, similar problems have been investigated in discrete settings, such as finite graphs and Markov chains. Typically, Chow, Huang, Li, and Zhou \cite{chow2012arma} investigated the relationships among three concepts defined on graphs: the free energy functional, the Fokker-Planck equation, and stochastic processes. It is well known that the notation of gradient flow makes sense only in context with an appropriate metric. As an alternative to the 2-Wasserstein metric defined on the continuous setting, several new metrics on the positive probability distributions with a finite graph as an underling space were constructed in \cite{chow2012arma}. Different choices for metric result in different Fokker-Planck equations. From the free energy viewpoint, they deduced a system of nonlinear ordinary differential equations, called the Fokker-Planck equation on graphs, which is the gradient flow of the free energy functional defined on a Riemannian manifold of positive probability distributions. From the stochastic viewpoint, they introduced a new interpretation of white noise perturbations to a Markov process on the discrete space, and derived another Fokker-Planck equation as the time evolution equation for its probability density function, which is not the same as the one obtained from the free energy functional. Under those settings, the unique global equilibrium of those Fokker-Planck equations are Gibbs distribution. Building on the framework constructed in \cite{chow2012arma}, the authors of \cite{che2016jde} proved the exponential rate of convergence towards the global equilibrium of these Fokker-Planck equations, measured both in the $L^2$ norm and in the (relative) entropy. Using this convergence result, they also proved that two Talagrand-type inequalities hold true based on two different metrics introduced in \cite{chow2012arma}. In \cite{maas2011jfa}, Maas constructed a metric similar to, but different from, the 2-Wasserstein metric, defined via a discrete variant of the Benamou-Brenier formula \cite{benamou2000nm}. Maas showed that with respect to this metric, the law of the continuous time Markov chain evolves as the gradient flow for the entropy defined on a finite set. Erbar and Maas \cite{erbar2012arma} introduced a new notion of Ricci curvature that applies to Markov chains on discrete spaces under the metric constructed in \cite{maas2011jfa}. For a more general free energy functional consisting of a Boltzmann entropy, a linear potential and a quadratic interaction energy, Chow, Li and Zhou \cite{chow2018dcds} deduced the Fokker-Planck equation on graphs as a gradient flow of this free energy functional. Their metric endowed on the positive probability distributions of the graph is similar to the one introduced by Maas \cite{maas2011jfa}. In this article \cite{chow2018dcds}, Chow, Li and Zhou also proved the so called Log-Sobolev inequality by using the convergence of the solution. The asymptotic properties of the solution also are studied. Several numerical examples related to similar topics are provided in \cite{chow2019jdde}.

All the work mentioned above was developed in the finite setting, such as finite graphs and finite Markov chains. To the best of the author's knowledge, no related results exist in the context of the infinite settings. There are two primary obstacles arise when attempting to extend the results from the finite case to more general infinite graphs. The first obstacle is how to construct the gradient flow on the infinite dimensional Riemannian manifold $(\mathcal{P}_0(G), \mathcal{W}_2)$. Unlike the case of finite graphs, where it is natural to derive the gradient flow equation from the inner product structure in the finite-dimensional Riemannian manifold. Establishing the gradient flow equation on infinite-dimensional manifolds requires more careful analysis and a deeper investigation into the properties of the individual components. The second obstacle is the existence and asymptotic behavior of the global solution to the Fokker-Planck equation, an infinite-dimensional ordinary differential equation, for all time $t>0$ in the positive probability density space. The proof of existence in finite graphs is elegant \cite{chow2012arma}, where the method involves constructing a carefully chosen bounded subset, which is compact and entirely contained within the probability space $\mathcal{P}_0(G)$, and demonstrating that the solution remains within this bounded subset if the initial data lie within it. This method has been applied in several studies \cite{chow2012arma, che2016jde, chow2018dcds, chow2019jdde}. Unfortunately, this approach is only applicable to finite graphs. In the infinite case, their definition of bounded sets become ambiguous at infinity. More importantly, bounded sets in infinite-dimensional spaces are not necessarily compact. The proof of asymptotic behavior for finite graphs used the traditional method for studying gradient flows in Euclidean space. This method is also not applicable in the case of infinite graphs.

In this paper, we are trying to overcome those obstacles. Before introducing the work presented in this paper, it is important to emphasize that all the results in this paper are based on the assumption that the weighted measure defined on the vertexes of the locally finite graph satisfies
\begin{align*}
\Pi:=\sum_{i=1}^\infty\pi(x_i)=\sum_{i=1}^\infty\pi_i=\infty,\quad\quad\alpha:=\inf_{i\in\mathbb{N}}\pi_i>0.
\end{align*}
The tangent space of any finite dimensional manifold is isomorphism to Euclidean space with the same dimension, whereas the tangent space of  an infinite dimensional manifold is isomorphism to some Banach or Hilbert space. This fact provides more flexibility in constructing the tangent space of the probability density space defined on the infinite graphs, allowing it to possess the necessary properties. We denote $\mathcal{P}_0(G)$ the positive probability density space, which is the interior of the probability density space $\mathcal{P}(G)$, one can find the detailed definition in Section \ref{sec:pre}. Define the following tangent space
\begin{align*}
T_{\bm{\rho}}\mathcal{P}_0(G)=\left\{\bm{\sigma}=(\sigma_i)_{i=1}^\infty\in\textrm{Ran}(B_{\bm{\rho}})\mid\sigma_i=\bm{\sigma}(x_i): V\to\mathbb{R},~\sum_{i=1}^\infty\sigma_i\pi_i=0\right\},
\end{align*}
where $B_{\bm{\rho}}$ is a negative weighted Laplacian operator and its definition will be given in Section \ref{sec:wasserstein}. Based on the definition of tangent space, the 2-Wasserstein distance was given by 
\begin{align*}
\mathcal{W}_2^2(\bm{\rho}^1,\bm{\rho}^2):=\inf\left\{\int_0^1\left\langle\dot{\bm{\rho}}, B_{\bm{\rho}}^{-1}\dot{\bm{\rho}}\right\rangle_{\bm{\pi}}\di t: \bm{\rho}(0)=\bm{\rho}^1, \bm{\rho}(1)=\bm{\rho}^2,\bm{\rho}\in\mathcal{C}\right\},
\end{align*}
where $\mathcal{C}$ is the set of all continuously differentiable curves $\bm{\rho}:[0,1]\to\mathcal{P}_0(G)$. The Fokker-Planck equation is given as 
\begin{align*}
\frac{\di \bm{\rho}}{\di t}=div_G(\bm{\rho}\nabla_G\bm{\Psi})+\Delta_G\bm{\rho}.
\end{align*}
Under the definition of the 2-Wasserstein metric, the Fokker-Planck equation can be viewed as a gradient flow of the free energy functional
\begin{align*}
\mathcal{F}(\bm{\rho})=\sum_{i=1}^\infty\pi_i\Psi_i\rho_i+\beta\sum_{i=1}^\infty\pi_i\rho_i\log{\rho_i},\quad\forall\bm{\rho}\in\mathcal{P}_0^*(G),
\end{align*}
on $(\mathcal{P}_0^*(G),\mathcal{W}_2)$, where $\bm{\Psi}=(\Psi_i)_{i=1}^\infty$ is the rapidly increasing potential of graph $G$ and $\beta>0$ is a universal constant. Without loss of generality, we assume $\beta=1$ in the following of this paper. The space $\mathcal{P}_0^*(G)$ is defined as a subset of the positive probability density space $\mathcal{P}_0(G)$ whose elements have finite second moment. This condition ensures sufficiently rapid decay of the probability densities, thereby guaranteeing that the free energy functional is well-defined. In mathematical terms, the space $\mathcal{P}_0^*(G)$ is defined as 
\begin{align*}
\mathcal{P}_0^*(G):=\left\{\bm{\rho}\in\mathcal{P}_0(G)\mid\sum_{i=1}^\infty d^2(x_1,x_i)\rho_i\pi_i<\infty,\quad\forall~x_i\in V\right\}.
\end{align*}
with some fixed $x_1\in V$,where $d(x_i,x_j)$ is the length of the shortest path form $x_i$ to $x_j$ on graph $G$.

To prove the global existence of the Fokker-Planck equation, we use the equation itself to derive a contradiction, showing that the solution will never touch the boundary of the probability density space. Then, we use the method of dividing the indicator set to find the decay properties of the solution to Fokker-Planck equation as a gradient flow to show that the solution asymptotically approaches the Gibbs density $\bm{\rho}*=(\rho_i^*)_{i=1}^\infty$, given by 
\begin{align}
\rho_i^*=\frac{1}{K}e^{-\Psi_i}~\textrm{with}~K=\sum_{j=1}^\infty \pi_je^{-\Psi_j},
\end{align}
which is in space $\mathcal{P}_0^*(G)$ because of the rapidly increase of the potential $\bm{\Psi}=(\Psi)_{i=1}^\infty$. Our main result is stated as follows:

\begin{theorem}\label{thm:globalexistence} 
Let \( G = (V, E, \bm{\pi}) \) be a connected, locally uniformly finite and stochastically complete graph. Then, the following statements hold:

(1) If $\bm{\rho}^0\in\mathcal{P}_0^*(G)$, the Fokker-Planck equation defined on the graph $G$ is a gradient flow on the infinite dimensional Riemannian manifold $(\mathcal{P}_0^*(G),\mathcal{W}_2)$;

(2) If $\bm{\rho}^0\in\mathcal{P}_0^*(G)$, the Fokker-Planck equation defined on the graph $G$ has a unique global solution in \(\mathcal{P}_0^*(G)\);

(3) For $r\in [2,\infty]$, the unique global solution $\bm{\rho}$ converges to the Gibbs density \(\bm{\rho}^* = (\rho_i^*)_{i=1}^\infty\in\mathcal{P}_0^*(G)\) under the $\ell^{r}(V,\bm{\pi})$ norm.
\end{theorem}

The rest of this paper is organized as follows. In Section \ref{sec:pre}, we introduce some basic settings about analysis on graphs. The construction of 2-Wasserstein metric and gradient flow equation on locally finite graphs will be presented in Section \ref{sec:wasserstein}. Finally, we study the Fokker-Plank equation on locally finite graphs, the global existence and asymptotic behavior of the solution will be proved in Section \ref{sec:FPE}.

\section{Preliminaries}\label{sec:pre}
Let $G=(V,E,\bm{\pi})$ be a connected, locally uniformly finite and stochastically complete weighted graph with vertex set $V$ and a fixed measure $\bm{\pi}=(\pi(x_i))_{i=1}^\infty=(\pi_i)_{i=1}^\infty$ on $V$ that satisfies
\begin{align*}
\Pi:=\sum_{i=1}^\infty\pi_i=\infty,\quad\quad\alpha:=\inf_{i\in\mathbb{N}}\pi_i>0,
\end{align*}
and the growth assumption, that is there exists a universal constant $K>0$ such that
\begin{align}\label{app:growth}
\pi_j\leq K\pi_i,
\end{align}
for any adjacent vertices $x_i$ and $x_j$ in $V$.
The set $E$ is the edge set, and $\bm{\omega}=(\omega_{ij})_{x_i,x_j\in V}$ contains the weight of each edge,
\begin{equation*}
\omega_{ij}=
\left\{\begin{aligned}
&\omega_{ij}>0,&&\textrm{if}~\{x_i,x_j\}\in E;&\\
&0,&&\textrm{otherwise}.&
\end{aligned}
\right.
\end{equation*}
The weight function satisfies $\omega_{ij}=\omega_{ji}$ for any $i,j\in\mathbb{N}$, and
\begin{align*}
\lambda:=\inf_{i,j\in\mathbb{N}\atop i\sim j}\omega_{ij}>0\quad\textrm{and}\quad\Lambda:=\sup_{i,j\in\mathbb{N}}\omega_{ij}<\infty,
\end{align*}
where $i\sim j$ implies $x_i$ and $x_j$ is adjacent for $i,j\in\mathbb{N}$.
We denote
\begin{align*}
N(i)=\{j\in \mathbb{N}\mid~\{x_i,x_j\}\in E\}
\end{align*}
as the set of vertices adjacent to $x_i$. We assume that for any $i\in\mathbb{N}$, $|N(i)|\leq C_V$ for some universal constant $C_V\geq 1$. The measure $\bm{\pi}=(\pi_i)_{i=1}^\infty$ and the edge weight function $\bm{\omega}$ satisfy the following relation:
\begin{align*}
\sum_{j\in \mathbb{N}}\omega_{ij}=\sum_{j\in N(i)}\omega_{ij}=\pi_i.
\end{align*}
We assume that the graph $G$ is an undirected graph with no loops or multiple edges. 

Let $C(V)$ be the set of real functions on $V$ and $C(V^2)$ be the set of real function on $V\times V$. For any $1\leq p<\infty$, we denote:
\begin{align*}
\ell^p(V,\bm{\pi}):=\left\{f\in C(V)\mid~\sum_{i\in\mathbb{N}}\pi_i\left|f(x_i)\right|^p<\infty\right\}
\end{align*}
as the set of $\ell^p$ integrable functions on $V$ with respect to the measure $\bm{\pi}$. For $p=\infty$,
\begin{align*}
\ell^{\infty}(V):=\left\{f\in C(V)\mid~\sup_{i\in\mathbb{N}}|f(x_i)|<\infty\right\}.
\end{align*}
The standard inner product is defined by
\begin{align}\label{def:innerproduct1}
\langle f,g\rangle_{\bm{\pi}}:=\sum_{i\in\mathbb{N}}f(x_i)g(x_i)\pi_i,\quad\forall f,g\in\ell^2(V,\bm{\pi}).
\end{align}
This makes $\ell^2(V,\bm{\pi})$ a Hilbert space.
Let $\phi\in C(V)$, define $\nabla_G\phi\in C(V^2)$ by
\begin{equation*}
\nabla_G\phi(x_i,x_j):=
\left\{\begin{aligned}
&\phi(x_j)-\phi(x_i),&&\textrm{if}~\{x_i,x_j\}\in E;&\\
&0,&&\textrm{otherwise}.&
\end{aligned}
\right.
\end{equation*}
The $\bm{\pi}$-Laplacian operator $\Delta_G$ on graph is defined by
\begin{align}\label{def:delta1}
\Delta_G\phi(x_i):=\sum_{j\in N(i)}\frac{\omega_{ij}}{\pi_i}\left(\phi(x_j)-\phi(x_i)\right),\quad\forall\phi\in C(V).
\end{align}
For any function $\Phi\in C(V^2)$, the divergence $div_G\Phi\in C(V)$ is defined by
\begin{align*}
div_G\Phi(x_i):=\sum_{j\in N(i)}\frac{\omega_{ij}}{\pi_i}\Phi(x_i,x_j).
\end{align*}
By the definition of gradient, divergence and Laplacian operator on $C(V)$, we have
\begin{align*}
div_G\left(\nabla_G\phi\right)=\Delta_G\phi,\quad\forall \phi\in C(V).
\end{align*}
Actually, direct calculation shows
\begin{equation*}
\begin{aligned}
div_G\left(\nabla_G\phi\right)(x_i)=&\sum_{j\in N(i)}\frac{\omega_{ij}}{\pi_i}\nabla_G\phi(x_i,x_j)\\
=&\sum_{j\in N(i)}\frac{\omega_{ij}}{\pi_i}\left(\phi(x_j)-\phi(x_i)\right)=\Delta_G\phi(x_i).
\end{aligned}
\end{equation*}

The probability density space on the graph $G$ is defined as follows
\begin{align*}
\mathcal{P}(G)=\left\{\bm{\rho}=(\rho_i)_{i=1}^\infty\mid~\rho_i=\rho(x_i): V\to\mathbb{R}, \sum_{i=1}^\infty\rho_i\pi_i=1~\textrm{and}~\rho_i\geq 0, \forall i\in\mathbb{N}\right\}.
\end{align*}
The boundary of $\mathcal{P}(G)$ is defined by
\begin{align*}
\partial\mathcal{P}(G):=\left\{\bm{\rho}=(\rho_i)_{i=1}^\infty\in\mathcal{P}(G)\mid\exists~i_0\in\mathbb{N},~\textrm{s.t.}~\rho_{i_0}=0\right\}.
\end{align*}
The positive probability density space (the interior of $\mathcal{P}(G)$) on $G$ is defined by $\mathcal{P}_0(G):=\mathcal{P}(G)/\partial\mathcal{P}(G)$.

For any probability density $\bm{\rho}\in\mathcal{P}(G)$ and functional $\Phi\in C(V^2)$, the weighted divergence is defined by
\begin{align*}
div_G(\bm{\rho}\Phi)(x_i):=\sum_{j\in N(i)}\frac{\omega_{ij}}{\pi_i}\Phi(x_i,x_j)\hat{\rho}(x_i,x_j),
\end{align*}
where $\hat{\rho}(x_i,x_j)$ is the Logarithmic mean
\begin{equation*}
\hat{\rho}(x_i,x_j):=
\left\{\begin{aligned}
&\frac{\rho(x_i)-\rho(x_j)}{\log\rho(x_i)-\log\rho(x_j)},&&\textrm{if}~\{x_i,x_j\}\in E,~\textrm{and}~\rho(x_i)> 0, \rho(x_j)> 0;&\\
&0,&&\textrm{others},&
\end{aligned}
\right.
\end{equation*}
for any $x_i,x_j\in V$. Notice that $\hat{\rho}$ is bounded for any $\bm{\rho}\in\mathcal{P}(G)$ under our basic setting. According to the definition of $\hat{\rho}$, the $\bm{\pi}$-Laplacian operator can be rewritten as
\begin{align*}
\Delta_G\rho(x_i)=\sum_{j\in N(i)}\frac{\omega_{ij}}{\pi_i}\left(\log\rho(x_j)-\log\rho(x_i)\right)\hat{\rho}(x_i,x_j),
\end{align*}
for $\bm{\rho}\in\mathcal{P}_0(G)$.

We denote $\bm{\Psi}=(\Psi(x_i))_{i=1}^\infty=(\Psi_i)_{i=1}^\infty$ as the potential on the graph $G$ (i.e., $\Psi_i$ is the potential at the state $x_i$). We assume that the potential such that
\begin{enumerate}
\item $\Psi_i\to\infty$ as $i\to\infty$;
\item The sequence is increasing rapidly enough as $i\to\infty$;
\item $|\nabla_G\bm{\Psi}|\leq C_{\Psi}$ with $C_{\Psi}>1$ is a constant.
\end{enumerate}
The conditions $(1)$ and $(2)$ are to ensure the Gibbs distribution lies in the space $\mathcal{P}_0^*(G)$, but the increase rate is not fast arbitrary because of condition $(3)$. The condition $(3)$ here is to guarantee that the Fokker-Planck equation possesses sufficiently desirable properties.

A function $h:(0,+\infty)\times V\times V\to\mathbb{R}$ is called a fundamental solution to the following heat equation on graph $G=(V,E,\bm{\pi})$,
\begin{align}\label{eq:heatequation}
u_t=\Delta_G u,
\end{align}
if for any bounded initial data $u^0\in C(V)$, the function
\begin{align*}
u(t,x_i)=\sum_{j\in\mathbb{N}}h(t,x_i,x_j)u^0(x_j)\pi_j,\quad\forall i\in\mathbb{N}, t>0,
\end{align*}
is differentiable in time variable $t$, satisfies the heat equation \eqref{eq:heatequation}, and for any $i\in\mathbb{N}$, there are $\lim_{t\to 0^+}u(t,x_i)=u^0(x_i)$. The heat semigroup associated with $-\Delta_G$ was given by 
\begin{align*}
e^{t\Delta_G}f(x_i):=\sum_{j\in\mathbb{N}}h(t,x_i,x_j)f(x_j)\pi_j,
\end{align*}
for all $i\in\mathbb{N}$ and $f\in C(V)$. It is  known that $\sum_{j\in\mathbb{N}}h(t,x_i,x_j)\pi(x_j)\leq 1$. The graph is called stochastically complete, if the following condition holds:
\begin{align*}
\sum_{j\in\mathbb{N}}h(t,x_i,x_j)\pi_j=1,\quad\forall t>0.
\end{align*}
Stochastically completeness ensures the conservation of mass within the graph, preventing any occurrence of mass leakage. For example, the solution of heat equation given by
\begin{align*}
u_i(t)=&\sum_{j\in\mathbb{N}}h(t,x_i,x_j)u^0(x_j)\pi_j.
\end{align*}
Hence, if the graph is stochastically incomplete, we will obtain
\begin{align*}
\sum_{i\in\mathbb{N}}\pi_iu_i(t)=\sum_{i\in\mathbb{N}}\pi_i\sum_{j\in\mathbb{N}}h(t,x_i,x_j)u^0(x_j)\pi_j=\sum_{j\in\mathbb{N}}u^0(x_j)\pi_j\sum_{i\in\mathbb{N}}\pi_ih(t,x_i,x_j)<1.
\end{align*}
This implies that the solution flows out of the probability density space \(\mathcal{P}(G)\), violating the conservation of total probability. One can find more details on this topic in \cite{horn2019jran, delmotte1999} and the references therein.

In the following of this paper, we will use both symbols $\bm{\rho}(x_i)$ and $\rho_i$ interchangeably, emphasizing here that these two notations are same. Similarly, if other quantities are written in analogous ways, they are also same, for example $\bm{\pi}(x_i)=\pi_i$ and $\bm{\Psi}(x_i)=\Psi_i$.

\section{Wasserstein type distance on locally uniformly finite graphs}\label{sec:wasserstein}
In this section, we construct a 2-Wasserstein type metric and derivative the gradient flow equation on the locally uniformly finite graph $G=(V,E,\bm{\pi})$. This section closely parallels the corresponding part in \cite{chow2012arma,chow2018dcds,maas2011jfa}, where the authors focused on finite graphs. Here, we deal with the infinite case.

\subsection{Wsserstein type distance I}
In this subsection, we construct the Wasserstein type distance from the continuous equation on locally finite graph. For any $\bm{\rho}\in\mathcal{P}(G)$, define the space 
\begin{align*}
L_{\bm{\rho}}^2(V^2,\bm{\pi}):=\left\{\Phi\in C(V^2)\mid\frac{1}{2}\sum_{i\in\mathbb{N}}\sum_{j\in N(i)}\frac{\omega_{ij}}{\pi_i}\Phi(x_i,x_j)^2\hat{\rho}(x_i,x_j)\pi_i<\infty\right\}.
\end{align*}
Given two vector fields $\Phi,\tilde{\Phi}\in L_{\bm{\rho}}^2(V^2,\bm{\pi})$ and $\bm{\rho}\in\mathcal{P}(G)$, the discrete weighted inner product is defined as follows
\begin{equation}\label{def:innerproduct2}
\begin{aligned}
\left(\Phi,\tilde{\Phi}\right)_{\bm{\rho}}:=&\sum_{\{x_i,x_j\}\in E}\frac{\omega_{ij}}{\pi_i}\Phi(x_i,x_j)\tilde{\Phi}(x_i,x_j)\hat{\rho}(x_i,x_j)\pi_i\\
=&\frac{1}{2}\sum_{i\in\mathbb{N}}\sum_{j\in N(i)}\frac{\omega_{ij}}{\pi_i}\Phi(x_i,x_j)\tilde{\Phi}(x_i,x_j)\hat{\rho}(x_i,x_j)\pi_i.
\end{aligned}
\end{equation}
In particular,
\begin{align*}
\left(\Phi,\Phi\right)_{\bm{\rho}}=\sum_{\{x_i,x_j\}\in E}\frac{\omega_{ij}}{\pi_i}\Phi(x_i,x_j)^2\hat{\rho}(x_i,x_j)\pi_i=\frac{1}{2}\sum_{i\in\mathbb{N}}\sum_{j\in N(i)}\frac{\omega_{ij}}{\pi_i}\Phi(x_i,x_j)^2\hat{\rho}(x_i,x_j)\pi_i.
\end{align*}
This inner product makes the space $L_{\bm{\rho}}^2(V^2,\bm{\pi})$ a Hilbert space.

For any $\bm{\rho}\in\mathcal{P}(G)$, we define
\begin{align*}
\mathcal{M}_{\bm{\rho}}(V^2,\bm{\pi}):=\left\{\Phi\in L_{\bm{\rho}}^2(V^2,\bm{\pi})\mid div_G(\bm{\rho}\Phi)\in\ell^1(V,\bm{\pi}),~\Phi~\textrm{is antisymmetry}\right\}.
\end{align*}
We also define the space
\begin{align*}
\mathcal{H}_{\bm{\rho}}(V,\bm{\pi}):=\left\{\bm{p}\in \ell^\infty(V)\mid\nabla_G\bm{p}\in \mathcal{M}_{\bm{\rho}}(V^2,\bm{\pi})\right\}.
\end{align*}
Notice that $\nabla_G\bm{p}$ is antisymmetry and satisfies
$div_G(\bm{\rho}\nabla_G\bm{p})\in\ell^1(V,\bm{\pi})$ naturally. The later is because the locally uniformly finite construction of the graph and the boundedness of $\bm{p}$. Both $\mathcal{H}_{\bm{\rho}}(V,\bm{\pi})$ and $\mathcal{M}_{\bm{\rho}}(V^2,\bm{\pi})$ are the key spaces in our process to construct the gradient flow equation.

Next, we introduce a lemma. The identity proved in this lemma is analogous to integration by parts formula in continuous spaces.

\begin{lemma}\label{lmm:identityinnerproduct}
Let $\Phi\in\mathcal{M}_{\bm{\rho}}(V^2,\bm{\pi}) $. Then, for any $\phi\in\mathcal{H}_{\bm{\rho}}(V,\bm{\pi})$, the following identity holds
\begin{align}\label{inq:innerprodct1}
-\left\langle div_G(\bm{\rho}\Phi),\phi\right\rangle_{\bm{\pi}}=\left(\Phi,\nabla_G\phi\right)_{\bm{\rho}}.
\end{align}
Furthermore, the following property holds
\begin{align}\label{inq:innerproduct2}
\sum_{i=1}^\infty\left(div_G\left(\bm{\rho}\Phi\right)\right)(x_i)\pi_i=0.
\end{align}
\end{lemma}
\begin{proof}
First, because $div_G(\bm{\rho}\Phi)\in\ell^1(V,\bm{\pi})$ and $\phi\in\ell^\infty(V)$, there is 
\begin{align*}
\left\langle div_G(\bm{\rho}\Phi),\phi\right\rangle_{\bm{\pi}}<\infty. 
\end{align*}
Meanwhile, since $\Phi,\nabla_G\phi\in L_{\bm{\rho}}^2(V^2,\bm{\pi})$, there is
\begin{align*}
\left(\Phi,\nabla_G\phi\right)_{\bm{\rho}}<\infty.
\end{align*}
Next, we prove identity \eqref{inq:innerprodct1}.
On the one hand, by the definition of discrete inner product \eqref{def:innerproduct2}, we have
\begin{align*}
\left(\Phi,\nabla_G\phi\right)_{\bm{\rho}}
=&\frac{1}{2}\sum_{i\in\mathbb{N}}\sum_{j\in N(i)}\frac{\omega_{ij}}{\pi_i}\left(\phi(x_j)-\phi(x_i)\right)\Phi(x_i,x_j)\hat{\rho}(x_i,x_j)\pi_i\\
=&\frac{1}{2}\sum_{i\in\mathbb{N}}\sum_{j\in N(i)}\omega_{ij}\left(\phi(x_j)-\phi(x_i)\right)\Phi(x_i,x_j)\hat{\rho}(x_i,x_j).
\end{align*}
On the other hand, by the definition of standard inner product \eqref{def:innerproduct1}, we have
\begin{equation}\label{id:inner1}
\begin{aligned}
&\left\langle div_G(\rho\Phi),\phi\right\rangle_{\bm{\pi}}\\
=&\sum_{i\in\mathbb{N}}\sum_{j\in N(i)}\frac{\omega_{ij}}{\pi_i}\Phi(x_i,x_j)\hat{\rho}(x_i,x_j)\phi(x_i)\pi_i\\
=&\frac{1}{2}\sum_{i\in\mathbb{N}}\sum_{j\in N(i)}\omega_{ij}\Phi(x_i,x_j)\hat{\rho}(x_i,x_j)\phi(x_i)+\frac{1}{2}\sum_{i\in\mathbb{N}}\sum_{j\in N(i)}\omega_{ij}\Phi(x_i,x_j)\hat{\rho}(x_i,x_j)\phi(x_i).
\end{aligned}
\end{equation}
Interchanging the symbol $i$ and $j$ in the first term above, we have
\begin{align*}
\frac{1}{2}\sum_{i\in\mathbb{N}}\sum_{j\in N(i)}\omega_{ij}\Phi(x_i,x_j)\hat{\rho}(x_i,x_j)\phi(x_i)
=\frac{1}{2}\sum_{j\in\mathbb{N}}\sum_{i\in N(j)}\omega_{ji}\Phi(x_j,x_i)\hat{\rho}(x_j,x_i)\phi(x_j).
\end{align*}
By the symmetry of $\omega_{ij}=\omega_{ji}$, $\hat{\rho}(x_i,x_j)=\hat{\rho}(x_j,x_i)$, and the antisymmetry of $\Phi(x_i,x_j)=-\Phi(x_j,x_i)$, we obtain
\begin{align*}
\frac{1}{2}\sum_{i\in\mathbb{N}}\sum_{j\in N(i)}\omega_{ij}\Phi(x_i,x_j)\hat{\rho}(x_i,x_j)\phi(x_i)
=&-\frac{1}{2}\sum_{j\in\mathbb{N}}\sum_{i\in N(j)}\omega_{ij}\Phi(x_i,x_j)\hat{\rho}(x_i,x_j)\phi(x_j).
\end{align*}
Since $div_G(\bm{\rho}\Phi)\in\ell^1(V,\bm{\pi})$ and $\phi\in\ell^\infty(V)$, the series $\sum_{j\in\mathbb{N}}\sum_{i\in N(j)}\omega_{ij}\Phi(x_i,x_j)\hat{\rho}(x_i,x_j)\phi(x_j)$ is absolutely convergent. Hence, we can change the order of summation, that is
\begin{align*}
\frac{1}{2}\sum_{i\in\mathbb{N}}\sum_{j\in N(i)}\omega_{ij}\Phi(x_i,x_j)\hat{\rho}(x_i,x_j)\phi(x_i)
=&-\frac{1}{2}\sum_{j\in\mathbb{N}}\sum_{i\in N(j)}\omega_{ij}\Phi(x_i,x_j)\hat{\rho}(x_i,x_j)\phi(x_j)\\
=&-\frac{1}{2}\sum_{i\in\mathbb{N}}\sum_{j\in N(i)}\omega_{ij}\Phi(x_i,x_j)\hat{\rho}(x_i,x_j)\phi(x_j).
\end{align*}
Take this into the equality \eqref{id:inner1}, we have
\begin{align*}
&\left\langle div_G(\rho\Phi),\phi\right\rangle_{\bm{\pi}}\\
=&-\frac{1}{2}\sum_{i\in\mathbb{N}}\sum_{j\in N(i)}\omega_{ij}\Phi(x_i,x_j)\hat{\rho}(x_i,x_j)\phi(x_j)+\frac{1}{2}\sum_{i\in\mathbb{N}}\sum_{j\in N(i)}\omega_{ij}\Phi(x_i,x_j)\hat{\rho}(x_i,x_j)\phi(x_i)\\
=&-\frac{1}{2}\sum_{i\in\mathbb{N}}\sum_{j\in N(i)}\omega_{ij}\Phi(x_i,x_j)\hat{\rho}(x_i,x_j)\left(\phi(x_j)-\phi(x_i)\right).
\end{align*}
This implies the equivalence between $-\left\langle div_G(\rho\Phi),\phi\right\rangle_{\bm{\pi}}$ and $\left(\Phi,\nabla_G\phi\right)_{\bm{\rho}}$.
	
Identity \eqref{inq:innerproduct2} is a direct consequence of \eqref{inq:innerprodct1} by taking  $\phi=\bm{1}=(1,1,\cdots)\in\mathcal{H}_{\rho}(V,\bm{\pi})$.
\end{proof}

If we denote a new operator $\nabla_G^*: \mathcal{M}_{\bm{\rho}}(V^2,\bm{\pi})\to \ell^1(V,\bm{\pi})$ by $\nabla_G^*\Phi=-div_G(\bm{\rho}\Phi)$, and observe Lemma \ref{lmm:identityinnerproduct} from another perspective, we will find that
\begin{align*}
\left\langle \nabla_G^*\Phi,\phi\right\rangle_{\bm{\pi}}=\left(\Phi,\nabla_G\phi\right)_{\bm{\rho}}.
\end{align*}
This shows the gradient operator $\nabla_G$ is the adjoint operator of $\nabla_G^*$. Hence, we have
\begin{align}\label{id:Hodge}
\mathcal{M}_{\bm{\rho}}(V^2,\bm{\pi})=\textrm{Ran}(\nabla_G)\oplus^\perp\textrm{Ker}(\nabla_G^*).
\end{align}
As a consequence, we have the following lemma, which can be viewed as a discrete type Hodge's decomposition.

\begin{lemma}\label{lmm:hodge}
Given an function $\bm{v}\in \mathcal{M}_{\bm{\rho}}(V^2,\bm{\pi})$ on the graph $G$, and a probability density $\bm{\rho}\in\mathcal{P}_0(G)$, there exists a unique decomposition, such that 
\begin{align*}
\bm{v}=\nabla_G\bm{p}+\bm{u},\quad\textrm{and}\quad div_G\left(\bm{\rho} \bm{u}\right)=0,
\end{align*}
where $\bm{p}\in\mathcal{H}_{\bm{\rho}}(V,\bm{\pi})$. In addition, the following property holds,
\begin{align*}
(\bm{v},\bm{v})_{\bm{\rho}}=\left(\nabla_G\bm{p},\nabla_G\bm{p}\right)_{\bm{\rho}}+(\bm{u},\bm{u})_{\bm{\rho}}.
\end{align*}
\end{lemma}

Next, we define the discrete analogue of 2-Wasserstein metric on the positive probability density space $\mathcal{P}_0(G)$ using the continuity equation. See \cite{benamou2000nm} for the original work in the continuous setting. For any $\bm{\rho}^1,\bm{\rho}^2\in\mathcal{P}_0(G)$, define
\begin{align*}
\mathcal{W}_1^2(\bm{\rho}^1,\bm{\rho}^2):=\inf\left\{\int_0^1(\bm{v}(t),\bm{v}(t))_{\rho(t)}\di t: \dot{\bm{\rho}}+div_G(\bm{\rho}(t)\bm{v}(t))=0, \bm{\rho}(0)=\bm{\rho}^1, \bm{\rho}(1)=\bm{\rho}^2\right\},
\end{align*}
where $\dot{\bm{\rho}}=\frac{\di}{\di t}\bm{\rho}$, and the infimum is taken over all antisymmetric function $\bm{v}\in \mathcal{M}_{\bm{\rho}}(V^2,\bm{\pi})$, and the continuously differentiable curve $\bm{\rho}:[0,1]\to\mathcal{P}_0(G)$. As a consequence of Lemma \ref{lmm:hodge}, the metric above can be rewritten as
\begin{equation}\label{def:w1}
\begin{aligned}
&\mathcal{W}_1^2(\bm{\rho}^1,\bm{\rho}^2):=\inf\left\{\int_0^1(\nabla_G\bm{p}(t),\nabla_G\bm{p}(t))_{\rho(t)}\di t: \right.\\
&\qquad\qquad\qquad\qquad\qquad\qquad\qquad\qquad\left.\dot{\bm{\rho}}+div_G(\bm{\rho}(t)\nabla_G\bm{p}(t))=0, \bm{\rho}(0)=\bm{\rho}^1, \bm{\rho}(1)=\bm{\rho}^2\right\},
\end{aligned}
\end{equation}
where the infimum is taken over all function $\bm{p}\in\mathcal{H}_{\bm{\rho}}(V,\bm{\pi})$ and the continuously differentiable curve $\bm{\rho}:[0,1]\to\mathcal{P}_0(G)$.

We define a new operator $A_{\bm{\rho}}$ from $\mathcal{H}_{\bm{\rho}}(V,\bm{\pi})$ to $C(V)$ by
\begin{align*}
A_{\bm{\rho}}\bm{p}(x_i):=-div_G(\rho\nabla_G\bm{p})(x_i)=-\sum_{j\in N(i)}\frac{\omega_{ij}}{\pi_i}\left(\bm{p}(x_j)-\bm{p}(x_i)\right)\hat{\rho}(x_i,x_j),
\end{align*}
for any $i\in\mathbb{N}$. Because the graph is locally finite, the operator $A_{\bm{\rho}}$ is well-defined.
Obviously, it is a negative weighted Laplacian operator, distinct from $-\Delta_G$. We denote by $A_{\bm{\rho}}^{-1}$ the pseudo-inverse operator of the weighted Laplacian operator $A_{\bm{\rho}}$. Actually, according to the definition of the operator $A_{\bm{\rho}}$, the constant function $\bm{1}=(1)_{i=1}^\infty$ is the eigenfunction of it with eigenvalue $0$. Note that
\begin{align*}
\left\langle A_{\bm{\rho}}\bm{p},\bm{p}\right\rangle_{\bm{\pi}}=\frac{1}{2}\sum_{i\in\mathbb{N}}\sum_{j\in N(i)}\frac{\omega_{ij}}{\pi_i}\left(\bm{p}(x_i)-\bm{p}(x_j)\right)^2\hat{\rho}(x_i,x_j)\pi_i=0,
\end{align*}
which indicates $\bm{p}(x_i)=\bm{p}(x_j)$ for $j\in N(i)$ and $i\in\mathbb{N}$. Since the graph $G$ is connected, it follows that $\bm{p}(x_i)=\bm{p}(x_j)$ for any $i,j\in\mathbb{N}$. Hence, $0$ is a simple eigenvalue, and $\textrm{Ker}(A_{\bm{\rho}})=\{\bm{c}:=(c)_{i=1}^\infty\mid c\in\mathbb{R}\}$. Let $\mathcal{R}$ be the quotient space $D(A_{\bm{\rho}})/\textrm{Ker}(A_{\bm{\rho}})=\mathcal{H}_{\bm{\rho}}(V,\bm{\pi})/\textrm{Ker}(A_{\bm{\rho}})$. In other words, for $\bm{p}\in C(V)$, we consider the equivalence class
\begin{align*}
[\bm{p}]=\left\{\left(p_1+c,p_2+c,\cdots,p_n+c,\cdots\right):~c\in\mathbb{R}\right\},
\end{align*}
and all such equivalent classes form the infinite dimensional space $\mathcal{R}$. We define the operator $\tilde{A}_{\bm{\rho}}:\mathcal{R}\to C(V)$ as follows
\begin{align*}
\tilde{A}_{\bm{\rho}}([\bm{p}])=A_{\bm{\rho}}(\bm{p}),\quad\forall \bm{p}\in \mathcal{K}_{\bm{\rho}}(V,\bm{\pi}).
\end{align*}
Obviously, the operator $\tilde{A}_{\bm{\rho}}$ is invertible. We denote its inverse operator by $\tilde{A}_{\bm{\rho}}^{-1}$ from $\textrm{Ran}(\tilde{A}_{\bm{\rho}})$ to $\mathcal{R}$. Furthermore, we have 
\begin{align*}
\mathcal{R}\times\textrm{Ker}(A_{\bm{\rho}})=\mathcal{H}_{\bm{\rho}}(V,\bm{\pi})/\textrm{Ker}(A_{\bm{\rho}})\times\textrm{Ker}(A_{\bm{\rho}})\cong \mathcal{H}_{\bm{\rho}}(V,\bm{\pi}).
\end{align*}
Hence, we define 
\begin{align*}
A_{\bm{\rho}}^{-1}(\bm{\sigma})=\bm{p},
\end{align*}
with $\tilde{A}_{\bm{\rho}}^{-1}(\bm{\sigma})=[\bm{p}]$. Here $\bm{p}$ is a representation of $[\bm{p}]$.

The tangent space of $\mathcal{P}_0(G)$ at $\bm{\rho}\in\mathcal{P}_0(G)$ was defined by
\begin{align*}
T_{\bm{\rho}}\mathcal{P}_0(G)=\left\{\bm{\sigma}=(\sigma_i)_{i=1}^\infty\in \textrm{Ran}(A_{\bm{\rho}})\mid\sigma_i=\bm{\sigma}(x_i): V\to\mathbb{R},~\sum_{i=1}^\infty\sigma_i\pi_i=0\right\}.
\end{align*}
Next, we present the equivalence between the tangent space $T_{\bm{\rho}}\mathcal{P}_0(G)$ and the range of operator $A_{\bm{\rho}}$ throughout the following lemma.

\begin{lemma}\label{lmm:operatorA}
For a given $\bm{\sigma}\in T_{\bm{\rho}}\mathcal{P}_0(G)$, there exists a unique real function $\bm{p}=(p_i)_{i=1}^\infty\in \mathcal{H}_{\bm{\rho}}(V,\bm{\pi})$, up to a constant shift, such that
\begin{align*}
\bm{\sigma}=A_{\bm{\rho}}\bm{p}=-div_G(\rho\nabla_G\bm{p}).
\end{align*}
Moreover, we have $T_{\bm{\rho}}\mathcal{P}_0(G)=\textrm{Ran}A_{\bm{\rho}}$.
\end{lemma}
\begin{proof}
On the one hand, for any $\bm{p}\in\mathcal{H}_{\bm{\rho}}(V,\bm{\pi})$, by \eqref{inq:innerproduct2} in Lemma \ref{lmm:identityinnerproduct}, we have
\begin{align*}
\sum_{i\in\mathbb{N}}A_{\bm{\rho}}\bm{p}(x_i)\pi_i=-\sum_{i=1}^\infty div_G\left(\rho\nabla_G\bm{p}\right)(x_i)\pi_i=0.
\end{align*}
This implies that $\textrm{Ran}(A_{\bm{\rho}})\subset T_{\bm{\rho}}\mathcal{P}_0(G)$.
	
On the other hand, the definition of $T_{\bm{\rho}}\mathcal{P}_0(G)$ yields $T_{\bm{\rho}}\mathcal{P}_0(G)\subset\textrm{Ran}(A_{\bm{\rho}})$. Hence, there is $T_{\bm{\rho}}\mathcal{P}_0(G)=\textrm{Ran}(A_{\bm{\rho}})$.
\end{proof}

\begin{remark}
In the finite graph case, i.e. $|V|=n<\infty$, the range $\textrm{Ran}(A_{\bm{\rho}})$ of operator $A_{\bm{\rho}}$ is a finite dimensional space, which is naturally a closed space. Hence, the tangent space is 
\begin{align*}
T_{\bm{\rho}}\mathcal{P}_0(G)=\left\{\bm{\sigma}=(\sigma_i)_{i=1}^\infty\in\mathbb{R}^n\mid\sigma_i=\bm{\sigma}(x_i): V\to\mathbb{R},~\sum_{i=1}^\infty\sigma_i\pi_i=0\right\},
\end{align*}
and $\textrm{Ker}(A_{\bm{\rho}})^\perp=\textrm{Ran}(A_{\bm{\rho}})=T_{\bm{\rho}}\mathcal{P}_0(G)$. One can find the details in \cite{chow2018dcds}.
\end{remark}

Next, we define an innerproduct on the tangent space $T_{\bm{\rho}}\mathcal{P}_0(G)$:

\begin{definition}\label{def:innerproduct}
Let $\bm{\sigma}^1, \bm{\sigma}^2\in T_{\bm{\rho}}\mathcal{P}_0(G)$, define the inner product $g_{\bm{\rho}}^{(1)}:T_{\rho}\mathcal{P}_0(G)\times T_{\rho}\mathcal{P}_0(G)\to\mathbb{R}$ by
\begin{align*}
g_{\bm{\rho}}^{(1)}(\bm{\sigma}^1,\bm{\sigma}^2):=\left\langle\bm{\sigma}^1, A_{\bm{\rho}}^{-1}\bm{\sigma}^2\right\rangle_{\bm{\pi}}=\left\langle A_{\bm{\rho}}\bm{p}^1,\bm{p}^2\right\rangle_{\bm{\pi}}=\left\langle\nabla_G\bm{p}^1,\nabla_G\bm{p}^2\right\rangle_{\bm{\rho}},
\end{align*}
where $\bm{\sigma}^1=A_{\bm{\rho}}\bm{p}^1$ and $\bm{\sigma}^2=A_{\bm{\rho}}\bm{p}^2$.
\end{definition}

As a consequence of Definition \ref{def:innerproduct}, the 2-Wasserstein distance in \eqref{def:w1} can be rewritten as
\begin{align*}
\mathcal{W}_1^2(\bm{\rho}^1,\bm{\rho}^2):=\inf\left\{\int_0^1g_{\bm{\rho}}^{(1)}(\dot{\bm{\rho}},\dot{\bm{\rho}})\di t: \bm{\rho}(0)=\bm{\rho}^1, \bm{\rho}(1)=\bm{\rho}^2,\bm{\rho}\in\mathcal{C}\right\},
\end{align*}
where $\mathcal{C}$ is the set of all continuously differentiable curves $\bm{\rho}(t):[0,1]\to\mathcal{P}_0(G)$.

Now, we give the gradient flow structure of any functional $\mathcal{J}$ from $\mathcal{P}_0(G)$ to $\mathbb{R}$ on the infinite dimensional Riemannian manifold $(\mathcal{P}_0(G),\mathcal{W}_1)$. 

\begin{lemma}\label{lmm:gradientflow}
For the connected, locally uniformly finite and stochastically complete graph $G$, the gradient flow of functional $\mathcal{J}\in C^2(\mathcal{P}(G))$ with $\frac{\delta}{\delta\bm{\rho}}\mathcal{J}\in\mathcal{H}_{\bm{\rho}}(V,\bm{\pi})$ on $(\mathcal{P}_0(G),\mathcal{W}_1)$ is
\begin{align}\label{eq:gradientflow1}
\dot{\bm{\rho}}(t)=-A_{\bm{\rho}}\frac{\delta}{\delta\bm{\rho}}\mathcal{J}(\bm{\rho}),
\end{align}
or equivalently,
\begin{align}\label{eq:gradientflow2}
\frac{\di\bm{\rho}}{\di t}=div_G\left(\rho\nabla_G\left(\frac{\delta}{\delta\bm{\rho}}\mathcal{J}(\bm{\rho})\right)\right).
\end{align}
\end{lemma}
\begin{proof}
For any $\bm{\sigma}\in T_{\bm{\rho}}\mathcal{P}_0(G)$, there exists a unique real function $\bm{p}\in\mathcal{H}_{\bm{\rho}}(V,\bm{\pi})$, such that $\bm{\sigma}=A_{\bm{\rho}}\bm{p}=-div_G(\bm{\rho}\nabla_G\bm{p})$. On the one hand, by the definition of the inner product $g_{\bm{\rho}}^{(1)}(\cdot,\cdot)$,  we have
\begin{align*}
g_{\bm{\rho}}^{(1)}(\dot{\bm{\rho}},\bm{\sigma})=\left\langle\dot{\bm{\rho}},\bm{p}\right\rangle_{\bm{\pi}}.
\end{align*}
On the other hand,
\begin{align*}
\left\langle\frac{\delta}{\delta\bm{\rho}} \mathcal{J},\bm{\sigma}\right\rangle_{\bm{\pi}}=\left\langle A_{\bm{\rho}}\frac{\delta}{\delta\bm{\rho}} \mathcal{J},\bm{p}\right\rangle_{\bm{\pi}}.
\end{align*}
According to the definition of gradient flow on a Riemannian manifold
\begin{align*}
g_{\bm{\rho}}(\dot{\bm{\rho}},\bm{\sigma})=-\left\langle\frac{\delta}{\delta\bm{\rho}} \mathcal{J},\bm{\sigma}\right\rangle_{\bm{\pi}},\quad\forall\bm{\sigma}\in T_{\bm{\rho}}\mathcal{P}_0(G),
\end{align*} 
we obtain
\begin{align*}
\left\langle\dot{\bm{\rho}},\bm{p}\right\rangle_{\bm{\pi}}=-\left\langle A_{\bm{\rho}}\frac{\delta}{\delta\bm{\rho}} \mathcal{J},\bm{p}\right\rangle_{\bm{\pi}}.
\end{align*}
Because $\ell^2(V,\bm{\pi})\subset\mathcal{H}_{\bm{\rho}}(V,\bm{\pi})$ and $\bm{p}\in \mathcal{H}_{\bm{\rho}}(V,\bm{\pi})$ is arbitrary, we have
\begin{align*}
\dot{\bm{\rho}}=-A_{\bm{\rho}}\frac{\delta}{\delta\bm{\rho}} \mathcal{J}.
\end{align*}
\end{proof}

\subsection{Wasserstein type distance II}
In this subsection, we introduce a new 2-Wasserstein-type distance $\mathcal{W}_2(\cdot,\cdot)$, which is similar to but distinct from the previously one $\mathcal{W}_1(\cdot,\cdot)$. The new one enables a greater number of equations to possess a gradient flow structure on $(\mathcal{P}_0(G),\mathcal{W}_2)$. The new distance has no relation to continuous equation, we define it directly.

For any $\bm{\rho}\in\mathcal{P}(G)$, we define
\begin{align*}
\mathcal{K}_{\bm{\rho}}(V,\bm{\pi}):=\left\{\bm{p}\in C(V)\mid div_G(\bm{\rho}\nabla_G\bm{p})\in\ell^1(V,\bm{\pi})\right\}.
\end{align*}
It is obvious that $\mathcal{H}_{\bm{\rho}}(V,\bm{\pi})\subset\mathcal{K}_{\bm{\rho}}(V,\bm{\pi})$. We define a new operator $B_{\bm{\rho}}$ from $\mathcal{K}_{\bm{\rho}}(V,\bm{\pi})$ to $\ell^1(V,\bm{\pi})$ by
\begin{align*}
B_{\bm{\rho}}\bm{p}(x_i):=-div_G(\rho\nabla_G\bm{p})(x_i)=-\sum_{j\in N(i)}\frac{\omega_{ij}}{\pi_i}\left(\bm{p}(x_j)-\bm{p}(x_i)\right)\hat{\rho}(x_i,x_j),
\end{align*}
for any $i\in\mathbb{N}$. Notice that the operators $B_{\bm{\rho}}$ and $A_{\bm{\rho}}$ differ only in their domains of definition. We denote by $B_{\bm{\rho}}^{-1}$ the pseudo-inverse operator of the weighted Laplacian operator $B_{\bm{\rho}}$.

The tangent space of $\mathcal{P}_0(G)$ at $\bm{\rho}\in\mathcal{P}_0(G)$ was defined by
\begin{align*}
T_{\bm{\rho}}\mathcal{P}_0(G)=\left\{\bm{\sigma}=(\sigma_i)_{i=1}^\infty\in \textrm{Ran}(B_{\bm{\rho}})\mid\sigma_i=\bm{\sigma}(x_i): V\to\mathbb{R},~\sum_{i=1}^\infty\sigma_i\pi_i=0\right\}.
\end{align*}
The next lemma is parallels to \ref{lmm:operatorA}, in which we present the equivalence between the tangent space $T_{\bm{\rho}}\mathcal{P}_0(G)$ and the range of operator $B_{\bm{\rho}}$.

\begin{lemma}\label{lmm:operatorB}
For a given $\bm{\sigma}\in T_{\bm{\rho}}\mathcal{P}_0(G)$, there exists a unique real function $\bm{p}=(p_i)_{i=1}^\infty\in \mathcal{K}_{\bm{\rho}}(V,\bm{\pi})$, up to a constant shift, such that
\begin{align*}
\bm{\sigma}=B_{\bm{\rho}}\bm{p}=-div_G(\rho\nabla_G\bm{p}).
\end{align*}
Moreover, we have $T_{\bm{\rho}}\mathcal{P}_0(G)=\textrm{Ran}B_{\bm{\rho}}$.
\end{lemma}
\begin{proof}
On the one hand, for any $\bm{p}\in\mathcal{K}_{\bm{\rho}}(V,\bm{\pi})$, 
because $\left(A_{\bm{\rho}}\bm{p}(x_i)\right)_{i=1}^\infty\in\ell^1(V,\bm{\pi})$ and 
the antisymmetry of $\nabla_G\bm{p}$, 
\begin{align*}
\sum_{i\in\mathbb{N}}B_{\bm{\rho}}\bm{p}(x_i)\pi_i=-\sum_{i=1}^\infty div_G\left(\rho\nabla_G\bm{p}\right)(x_i)\pi_i=0.
\end{align*}
 This implies that $\textrm{Ran}(B_{\bm{\rho}})\subset T_{\bm{\rho}}\mathcal{P}_0(G)$.

On the other hand, the definition of $T_{\bm{\rho}}\mathcal{P}_0(G)$ yields $T_{\bm{\rho}}\mathcal{P}_0(G)\subset\textrm{Ran}(B_{\bm{\rho}})$. Hence, there is $T_{\bm{\rho}}\mathcal{P}_0(G)=\textrm{Ran}(B_{\bm{\rho}})$.
\end{proof}

Next, we define an innerproduct on the tangent space $T_{\bm{\rho}}\mathcal{P}_0(G)$:

\begin{definition}\label{def:innerproductB}
Let $\bm{\sigma}^1, \bm{\sigma}^2\in T_{\bm{\rho}}\mathcal{P}_0(G)$, define the inner product $g_{\bm{\rho}}^{(2)}:T_{\rho}\mathcal{P}_0(G)\times T_{\rho}\mathcal{P}_0(G)\to\mathbb{R}$ by
\begin{align*}
g_{\bm{\rho}}^{(2)}(\bm{\sigma}^1,\bm{\sigma}^2):=\left\langle\bm{\sigma}^1, B_{\bm{\rho}}^{-1}\bm{\sigma}^2\right\rangle_{\bm{\pi}}=\left\langle B_{\bm{\rho}}\bm{p}^1,\bm{p}^2\right\rangle_{\bm{\pi}}=\left\langle\nabla_G\bm{p}^1,\nabla_G\bm{p}^2\right\rangle_{\bm{\rho}},
\end{align*}
where $\bm{\sigma}^1=B_{\bm{\rho}}\bm{p}^1$ and $\bm{\sigma}^2=B_{\bm{\rho}}\bm{p}^2$.
\end{definition}

As a consequence of Definition \ref{def:innerproductB}, we define the 2-Wasserstein distance as
\begin{align*}
\mathcal{W}_2^2(\bm{\rho}^1,\bm{\rho}^2):=\inf\left\{\int_0^1g_{\bm{\rho}}^{(2)}(\dot{\bm{\rho}},\dot{\bm{\rho}})\di t: \bm{\rho}(0)=\bm{\rho}^1, \bm{\rho}(1)=\bm{\rho}^2,\bm{\rho}\in\mathcal{C}\right\},
\end{align*}
where $\mathcal{C}$ is the set of all continuously differentiable curves $\bm{\rho}(t):[0,1]\to\mathcal{P}_0(G)$.

Now, we give the gradient flow structure of any functional $\mathcal{J}$ from $\mathcal{P}_0(G)$ to $\mathbb{R}$ on the infinite dimensional Riemannian manifold $(\mathcal{P}_0(G),\mathcal{W}_2)$. 

\begin{lemma}\label{lmm:gradientflowB}
For the connected, locally uniformly finite and stochastically complete graph $G$, the gradient flow of functional $\mathcal{J}\in C^2(\mathcal{P}(G))$ with $\frac{\delta}{\delta\bm{\rho}}\mathcal{J}\in\mathcal{K}_{\bm{\rho}}(V,\bm{\pi})$ on $(\mathcal{P}_0(G),\mathcal{W}_2)$ is
\begin{align}\label{eq:gradientflowB1}
\dot{\bm{\rho}}(t)=-B_{\bm{\rho}}\frac{\delta}{\delta\bm{\rho}}\mathcal{J}(\bm{\rho}),
\end{align}
or equivalently,
\begin{align}\label{eq:gradientflowB2}
\frac{\di\bm{\rho}}{\di t}=div_G\left(\rho\nabla_G\left(\frac{\delta}{\delta\bm{\rho}}\mathcal{J}(\bm{\rho})\right)\right).
\end{align}
\end{lemma}
\begin{proof}
For any $\bm{\sigma}\in T_{\bm{\rho}}\mathcal{P}_0(G)$, there exists a unique real function $\bm{p}\in\mathcal{K}_{\bm{\rho}}(V,\bm{\pi})$, such that $\bm{\sigma}=B_{\bm{\rho}}\bm{p}=-div_G(\bm{\rho}\nabla_G\bm{p})$. On the one hand, by the definition of the inner $g_{\bm{\rho}}^{(2)}(\cdot,\cdot)$,  we have
\begin{align*}
g_{\bm{\rho}}^{(2)}(\dot{\bm{\rho}},\bm{\sigma})=\left\langle\dot{\bm{\rho}},\bm{p}\right\rangle_{\bm{\pi}}.
\end{align*}
On the other hand,
\begin{align*}
\left\langle\frac{\delta}{\delta\bm{\rho}} \mathcal{J},\bm{\sigma}\right\rangle_{\bm{\pi}}=\left\langle B_{\bm{\rho}}\frac{\delta}{\delta\bm{\rho}} \mathcal{J},\bm{p}\right\rangle_{\bm{\pi}}.
\end{align*}
According to the definition of gradient flow on a Riemannian manifold
\begin{align*}
g_{\bm{\rho}}(\dot{\bm{\rho}},\bm{\sigma})=-\left\langle\frac{\delta}{\delta\bm{\rho}} \mathcal{J},\bm{\sigma}\right\rangle_{\bm{\pi}},\quad\forall\bm{\sigma}\in T_{\bm{\rho}}\mathcal{P}_0(G),
\end{align*} 
we obtain
\begin{align*}
\left\langle\dot{\bm{\rho}},\bm{p}\right\rangle_{\bm{\pi}}=-\left\langle B_{\bm{\rho}}\frac{\delta}{\delta\bm{\rho}} \mathcal{J},\bm{p}\right\rangle_{\bm{\pi}}.
\end{align*}
Because $\ell^2(V,\bm{\pi})\subset\mathcal{K}_{\bm{\rho}}(V,\bm{\pi})$ and $\bm{p}\in \mathcal{K}_{\bm{\rho}}(V,\bm{\pi})$ is arbitrary, we have
\begin{align*}
\dot{\bm{\rho}}=-B_{\bm{\rho}}\frac{\delta}{\delta\bm{\rho}} \mathcal{J}.
\end{align*}
\end{proof}

\section{Fokker-Planck equation on locally uniformly finite graphs: \\
	proof of Theorem \ref{thm:globalexistence}}\label{sec:FPE}
In this section, we demonstrate the Fokker-Planck equation, given by
\begin{align}\label{eq:FP}
\frac{\di\bm{\rho}}{\di t}=div_G\left(\bm{\rho}\nabla_G\bm{\Psi}\right)+\Delta_G\bm{\rho},
\end{align}
is a gradient flow on the infinite dimensional Riemannian manifold $(\mathcal{P}_0(G),\mathcal{W}_2)$, and prove the global existence, uniqueness, and asymptotic behavior of the solution, as stated in Theorem \ref{thm:globalexistence}.

Recall that $\bm{\Psi}=(\Psi_i)_{i=1}^\infty=(\Psi(x_i))_{i=1}^\infty$ is the potential on $V$ satisfying $|\nabla_G\Psi(x_i,x_j)|\leq C_{\Psi}$ for all $i,j\in\mathbb{N}$. The free energy functional on $\mathcal{P}_0(G)$ is given by
\begin{align}\label{def:freeenergy}
\mathcal{F}(\bm{\rho})=\sum_{i=1}^\infty\pi_i\Psi_i\rho_i+\sum_{i=1}^\infty\pi_i\rho_i\log{\rho_i},\quad\forall\bm{\rho}\in\mathcal{P}_0(G),
\end{align}
with 
\begin{align*}
\frac{\delta}{\delta\bm{\rho}}\mathcal{F}(\bm{\rho})=\left(\Psi_1+1+\log\rho_1,\cdots,\Psi_i+1+\log\rho_i,\cdots\right).
\end{align*}
It is obvious that the free energy functional $\mathcal{F}$ over the entire space of strictly positive probability densities $\mathcal{P}_0(G)$ is not well-defined, as it may take infinite values. Such cases are not meaningful for the analysis of the Fokker–Planck equation. To guarantee that the free energy functional is finite, it is necessary to impose sufficient decay conditions on the probability density. Therefore, we restrict our attention to the subset of probability densities with bounded second moment, i.e., for a fixed $x_1$, we define
\begin{align*}
\mathcal{P}_0^*(G):=\left\{\bm{\rho}\in\mathcal{P}_0(G)\mid\sum_{i=1}^\infty d^2(x_1,x_i)\rho_i\pi_i<\infty,~\forall x_i\in V\right\},
\end{align*}
where $d(x_j,x_i)$ is the length of the shortest path form $x_i$ to $x_j$ on graph $G$.

The next lemma implies that the finite of second moment ensures the potential function $\mathcal{F}_1(\bm{\rho}):=\sum_{i=1}^\infty\pi_i\Psi_i\rho_i$ and the entropy functional $\mathcal{F}_2(\bm{\rho}):=\sum_{i=1}^\infty\pi_i\rho_i\log{\rho_i}$ are also finite.
\begin{lemma}
For any $\bm{\rho}\in\mathcal{P}_0^*(G)$, we have
\begin{enumerate}
\item[(1)] $\sum_{i=1}^\infty\pi_i\Psi_i\rho_i<\infty$.
\item[(2)] $\sum_{i=1}^\infty\pi_i\rho_i\log{\rho_i}<\infty$.
\end{enumerate}
\end{lemma}
\begin{proof}
\textbf{Step I.} Because the gradient of $\bm{\Psi}$ is bounded, i.e., 
\begin{align*}
|\nabla_G\bm{\Psi}(x_i,x_j)|=|\bm{\Psi}(x_j)-\bm{\Psi}(x_i)|<C_{\Psi},
\end{align*}
the functional $\bm{\Psi}$ is Lipschitz continuity. Fixed a vertex $x_1$, then for any vertex $x_i$, there is a path in the graph $G$,
\begin{align*}
x_1=x_{i_0}\to x_{i_1}\to\cdots\to x_{i_k}=x_i.
\end{align*}
Then, we have
\begin{equation*}
\begin{aligned}
\bm{\Psi}(x_i)=&\bm{\Psi}(x_1)+\sum_{m=0}^{k-1}\left(\bm{\Psi}(x_{i_{m+1}})-\bm{\Psi}(x_{i_{m}})\right)\\
\leq&\bm{\Psi}(x_1)+C_{\Psi}k.
\end{aligned}
\end{equation*}
This inequality holds for any path from $x_1$ to $x_i$. Hence, we obtain
\begin{align}\label{inq:psiestimate}
\bm{\Psi}(x_i)\leq\bm{\Psi}(x_1)+C_{\Psi}d(x_1,x_i).
\end{align}
By \eqref{inq:psiestimate}, we have
\begin{align*}
\sum_{i=1}^\infty\pi_i\rho_i\Psi_i\leq&\sum_{i=1}^\infty\pi_i\rho_i\left(\Psi_1+C_{\Psi}d(x_1,x_i)\right)\\
=&\Psi_1+C_{\Psi}\sum_{i=1}^\infty\pi_i\rho_id(x_1,x_i).
\end{align*}
According to the Cauchy-Schwarz inequality, we have
\begin{align*}
\sum_{i=1}^\infty\pi_i\rho_i d(x_1,x_i)\leq\left(\sum_{i=1}^\infty\pi_i\rho_i\right)^{\frac{1}{2}}\left(\sum_{i=1}^\infty\pi_i\rho_i d^2(x_1,x_i)\right)^{\frac{1}{2}}.
\end{align*}
Hence,
\begin{align*}
\sum_{i=1}^\infty\pi_i\rho_i d(x_1,x_i)<\infty.
\end{align*}
As a consequence, we obtain
\begin{align*}
\sum_{i=1}^\infty\pi_i\rho_i\Psi_i\leq \Psi_1+C_{\Psi}\sum_{i=1}^\infty\pi_i\rho_i d(x_1,x_i)<\infty.
\end{align*}

\textbf{Step II.} Let $\mu_i=\pi_i\rho_i$, then
\begin{align*}
\sum_{i=1}^\infty\mu_i=1,\quad\mu_i>0.
\end{align*}
Hence, the entropy can be rewritten as 
\begin{align*}
\sum_{i=1}^\infty\pi_i\rho_i\log\rho_i=&\sum_{i=1}^\infty\mu_i\log\rho_i
=\sum_{i=1}^\infty\mu_i\log\mu_i-\sum_{i=1}^\infty\mu_i\log\pi_i.
\end{align*}
Next, we prove $\sum_{i=1}^\infty\mu_i\log\mu_i<\infty$ and $\sum_{i=1}^\infty\mu_i\log\pi_i<\infty$ respectively. We first prove the term $\sum_{i=1}^\infty\mu_i\log\mu_i$ is finite. Because $\sum_{i=1}^\infty\mu_i d^2(x_1,x_i)<\infty$, we have $\lim_{i\to\infty}\mu_i d^2(x_1,x_i)=0$. This shows that there exists a constant $C_0>0$ such that $\mu_i\leq \frac{C_0}{d^2(x_1,x_i)}$. Then, we have
\begin{align*}
|\mu_i\log\mu_i|=\mu_i|\log\mu_i|\leq \mu_i\left(2\log d(x_1,x_i)+\log C_0\right).
\end{align*} 
Hence,
\begin{align*}
\sum_{i=1}^\infty|\mu_i\log\mu_i|\leq&\sum_{i=1}^\infty\mu_i\left(2\log d(x_1,x_i)+\log C_0\right)\\
=&2\sum_{i=1}^\infty\mu_i\log d(x_1,x_i)+\log C_0\\
\leq&2\sum_{i=1}^\infty\mu_id(x_1,x_i)+\log C_0.
\end{align*}
According to the Cauchy-Schwartz inequality, we have
\begin{align*}
\sum_{i=1}^\infty\mu_i d(x_1,x_i)\leq \left(\sum_{i=1}^\infty\mu_i\right)^{\frac{1}{2}}\left(\sum_{i=1}^\infty\mu_i d^2(x_1,x_2)\right)^{\frac{1}{2}}<\infty.
\end{align*}
Hence,
\begin{align*}
\sum_{i=1}^\infty|\mu_i\log\mu_i|<\infty.
\end{align*}
For the term $\sum_{i=1}^\infty\mu_i\log\pi_i$, by the growth assumption \eqref{app:growth} for any two vertices $x_i$ and $x_j$ adjacent, we have
\begin{align*}
\pi_i\leq K^{d(x_1,x_i)}\pi_1.
\end{align*}
This is 
\begin{align}\label{inq:pi}
\log\pi_i\leq d(x_1,x_i)\log K+\log\pi_1.
\end{align}
Then,
\begin{equation*}
\begin{aligned}
\sum_{i=1}^\infty|\mu_i\log\pi_i|\leq&\sum_{i=1}^\infty\mu_id(x_1,x_i)|\log K|+\sum_{i=1}^\infty\mu_i|\log\pi_1|\\
\leq&|\log K| \left(\sum_{i=1}^\infty\mu_i\right)^{\frac{1}{2}}\left(\sum_{i=1}^\infty\mu_i d^2(x_1,x_2)\right)^{\frac{1}{2}}+|\log\pi_1|<\infty.
\end{aligned}
\end{equation*}
Finally, we obtain
\begin{align*}
\sum_{i=1}^\infty\pi_i\rho_i\log\rho_i=\sum_{i=1}^\infty\mu_i\log\mu_i-\sum_{i=1}^\infty\mu_i\log\pi_i<\infty.
\end{align*}
\end{proof}

Direct calculation shows that the Gibbs density $\bm{\rho}*=(\rho_i^*)_{i=1}^\infty$ is given by 
\begin{align}\label{def:gibbs}
\rho_i^*=\frac{1}{K}e^{-\Psi_i}~\textrm{with}~K=\sum_{j=1}^\infty \pi_je^{-\Psi_j},
\end{align}
which is the only global minimizer of the free energy functional $\mathcal{F}$ in $\mathcal{P}_0^*(G)$. Notice that we assume the sequence of potential $(\Psi_i)_{i=1}^\infty$ is increase rapidly to ensure the Gibbs density $\bm{\rho}*=(\rho_i^*)_{i=1}^\infty$ is in $\mathcal{P}_0^*(G)$.

\begin{remark}
In this section, we choose the Riemannian manifold $(\mathcal{P}_0^*(G),\mathcal{W}_2)$ as the space we study, but the manifold $(\mathcal{P}_0^*(G),\mathcal{W}_1)$. Because the metric $\mathcal{W}_2(\cdot,\cdot)$ enables more Fokker-Planck equations to possess a gradient flow structure on $(\mathcal{P}_0^*(G),\mathcal{W}_2)$. In fact, If we want the Fokker-Planck equation to be a gradient flow on the Riemannian manifold $(\mathcal{P}_0^*(G),\mathcal{W}_1)$, it is necessary to satisfy both conditions simultaneously that
\begin{enumerate}
\item $\frac{\delta}{\delta\bm{\rho}}\mathcal{F}(\bm{\rho})\in\ell^\infty(V)\subset\mathcal{H}_{\bm{\rho}}(V,\bm{\pi})$, by Lemma \ref{lmm:gradientflow}.
\item The Gibbs distribution $\bm{\rho}^*=(\rho_i^*)_{i=1}^\infty\in\mathcal{P}_0^*(G)$, i.e., the series $\bm{\rho}^*=(\rho_i^*)_{i=1}^\infty=\left(\frac{1}{K}e^{-\Psi_i}\right)_{i=1}^\infty$ has a bounded second moment.
\end{enumerate}
The restrictions imposed by conditions (1) and (2) mean that only very few choices of $\bm{\Psi}=(\Psi)_{i=1}^\infty$ are admissible. As a result, only a limited class of Fokker-Planck equations can be studied as gradient flows. However, in order for the Fokker-Planck equation to possess a gradient flow structure on $(\mathcal{P}_0^*(G),\mathcal{W}_2)$, it is sufficient to require $\frac{\delta}{\delta\bm{\rho}}\mathcal{F}(\bm{\rho})\in C(V)$. For this purpose, we prefer to conduct our study on $(\mathcal{P}_0^*(G),\mathcal{W}_2)$.
\end{remark}

Next, we divide the proof of Theorem \ref{thm:globalexistence} into three propositions, which respectively establish the gradient flow structure of the Fokker-Planck equation, the global existence of its solutions, and their asymptotic behavior.

\begin{proposition}\label{pro:gradientflow}
Let \( G = (V, E, \bm{\pi}) \) be a connected, locally uniformly finite and stochastically complete graph. Then, the Fokker-Planck equation defined on the graph $G$ is a gradient flow on the infinite dimensional Riemannian manifold $(\mathcal{P}_0^*(G),\mathcal{W}_2)$. Moreover, the Fokker-Planck equation can be written as 
\begin{align}\label{eq:FPE}
\frac{\di\rho_i}{\di t}
=\sum_{j\in N(i)}\frac{\omega_{ij}}{\pi_i}\left[(\Psi(x_j)+\log\rho(x_j))-\left(\Psi(x_i)+\log\rho(x_i)\right)\right]\hat{\rho}(x_i,x_j),\quad\forall i\in\mathbb{N}.
\end{align}
\end{proposition}
\begin{proof}
On the one hand, the Fokker-Planck equation is given by 
\begin{align}\label{eq:FPEoriginal}
\frac{\di\bm{\rho}}{\di t}=div_G\left(\bm{\rho}\nabla_G\bm{\Psi}\right)+\Delta_G\bm{\rho}.
\end{align}
A direct calculation shows that 
\begin{equation}\label{id:linearterm}
\begin{aligned}
div_G\left(\bm{\rho}\nabla_G\bm{\Psi}\right)(x_i)
=\sum_{j\in N(i)}\frac{\omega_{ij}}{\pi_i}\left[\Psi(x_j)-\Psi(x_i)\right]\hat{\rho}(x_i,x_j),\quad\forall i\in\mathbb{N},
\end{aligned}
\end{equation}
and 
\begin{equation}\label{id:laplaceterm}
\begin{aligned}
\Delta_G\bm{\rho}(x_i)
=\sum_{j\in N(i)}\frac{\omega_{ij}}{\pi_i}\left[(\log\rho(x_j)-\log\rho(x_i))\right]\hat{\rho}(x_i,x_j),\quad\forall i\in\mathbb{N}.
\end{aligned}
\end{equation}
Submitting \eqref{id:linearterm} and \eqref{id:laplaceterm} into the Fokker-Planck equation \eqref{eq:FPEoriginal}, we obtain
\begin{align}\label{eq:FPE1}
\frac{\di\rho_i}{\di t}
=\sum_{j\in N(i)}\frac{\omega_{ij}}{\pi_i}\left[(\Psi(x_j)+\log\rho(x_j))-\left(\Psi(x_i)+\log\rho(x_i)\right)\right]\hat{\rho}(x_i,x_j),\quad\forall i\in\mathbb{N}.
\end{align}

On the other hand, because $|N(i)|\leq C_V$ for any $i\in\mathbb{N}$ and the definition and boundedness of $\hat{\rho}$, it is easy to verify $\frac{\delta}{\delta\bm{\rho}}\mathcal{F}\in\mathcal{K}_{\bm{\rho}}(V,\bm{\pi})$. Taking $\mathcal{F}$ into the gradient flow equation \eqref{eq:gradientflow2}, we will obtain
\begin{align*}
\frac{\di\rho_i}{\di t}=\sum_{j\in N(i)}\frac{\omega_{ij}}{\pi_i}\left[(\Psi_j+\log\rho_j)-(\Psi_i+\log\rho_i)\right]\hat{\rho}(x_i,x_j), \quad\forall i\in\mathbb{N},
\end{align*}
which is same to \eqref{eq:FPE}.
\end{proof}

Next, we prove the global existence of the solution to Fokker-Planck equation in $\mathcal{P}_0^*(G)$.

\begin{proposition}\label{pro:existence}
Let \( G = (V, E, \bm{\pi}) \) be a connected, locally uniformly finite and stochastically complete graph. Then, for any initial data $\bm{\rho}^0 \in \mathcal{P}_0^*(G)$, the Fokker-Planck equation defined on the graph $G$ has a unique global solution in \(\mathcal{P}_0^*(G)\).
\end{proposition}
\begin{proof}
Firstly, by the boundedness of gradient of the potential, the definition of $\hat{\rho}$, and the locally uniformly finite structure of the graph, it is easy to observe that the following series is absolutely convergent for any $\bm{\rho}\in\mathcal{P}_0^*(G)$,
\begin{align*}
\sum_{i=1}^\infty\pi_i\frac{\di\rho_i}{\di t}=\sum_{i=1}^\infty\sum_{j\in N(i)}\omega_{ij}\left[(\Psi(x_j)+\log\rho(x_j))-(\Psi(x_i)+\log\rho(x_i))\right]\hat{\rho}(x_i,x_j).
\end{align*}
As a consequence, there holds
\begin{align*}
\sum_{i=1}^\infty\pi_i\frac{\di\rho_i}{\di t}=0,\quad\forall\bm{\rho}=(\rho_i)_{i=1}^\infty\in\mathcal{P}_0^*(G).
\end{align*}
Furthermore, by the gradient flow equation $\eqref{eq:gradientflow1}$, we have
\begin{align*}
\frac{\di}{\di t}\mathcal{F}(\bm{\rho})=&\left\langle\frac{\delta}{\delta\bm{\rho}}\mathcal{F}(\bm{\rho}),\frac{\di}{\di t}\bm{\rho}\right\rangle_{\bm{\pi}}\\
=&-\left\langle\frac{\delta}{\delta\bm{\rho}}\mathcal{F}(\bm{\rho}),B_{\bm{\rho}}\frac{\delta}{\delta\bm{\rho}}\mathcal{F}(\bm{\rho})\right\rangle_{\bm{\pi}}\\
=&-g_{\bm{\rho}}\left(\frac{\delta}{\delta\bm{\rho}}\mathcal{F}(\bm{\rho}),\frac{\delta}{\delta\bm{\rho}}\mathcal{F}(\bm{\rho})\right)<0.
\end{align*}
This implies the free energy is decrease along the solution of the Fokker-Planck equation. Hence, by Picard theorem, there exists a unique local solution $\bm{\rho}\in C^1([0,T_0); \mathcal{P}_0^*(G))$, if $\bm{\rho}_0\in\mathcal{P}_0^*(G)$.

Next, we prove the solution $\bm{\rho}(t)$ will never reach on the boundary $\partial\mathcal{P}(G)$ at any time $0\leq T<\infty$ by derive a contradiction. This yields the solution exists globally in the space $\mathcal{P}^*_0(G)$.

Assume the solution $\bm{\rho}$ hit the boundary at some time $0<T<\infty$ such that $T\geq T_0$ first time at the point $\bm{\mu}=(\mu_i)_{i=1}^\infty\in\partial\mathcal{P}(G)$ with 
\begin{align*}
\mu_i=0,\quad\textrm{for}~i\in M_1,
\end{align*}
where $M_1$ is the largest subset of $\mathbb{N}$ that ensures $\mu_i=0$ for $i\in M_1$. We denote $M_2=\mathbb{N}/M_1$.
Due to the connectivity of the graph $G$, there exists at least one index $i\in M_1$ and one index $k\in M_2$ such that $x_i$ and $x_k$ are adjacent. Hence, we have
\begin{align}\label{eq:sum}
\lim_{t\to T^-}\left(\log\rho_k(t)-\log\rho_i(t)\right)\hat{\rho}(x_i,x_k)=\lim_{t\to T^-}\left(\rho_k(t)-\rho_i(t)\right)>0,
\end{align}
and 
\begin{align*}
\lim_{t\to T^-}(\Psi_k-\Psi_i)\hat{\rho}(x_i,x_k)=0.
\end{align*}
Hence, we have
\begin{align*}
\left[(\Psi_k+\log\rho_k(T))-(\Psi_i+\log\rho_i(T))\right]\hat{\rho}(x_i,x_k)(T)>0.
\end{align*}
For the case $i,j\in M_1$, we have 
\begin{align*}
\lim_{t\to T^-}(\rho(x_j)-\rho(x_i))=0,\quad\textrm{and}\quad\lim_{t\to T^-}\hat{\rho}(x_i,x_j)=0.
\end{align*}
Hence, we have
\begin{align*}
&\left[(\Psi_j+\log\rho_j(T))-(\Psi_i+\log\rho_i(T))\right]\hat{\rho}(x_i,x_j)(T)\\
&\qquad\qquad\qquad\qquad\qquad\qquad=\left[(\Psi_j-\Psi_i)\hat{\rho}(x_i,x_j)(T)+(\rho_j(T)-\rho_i(T))\right]=0.
\end{align*}
Taking these into the Fokker-Planck equation, we obtain
\begin{align*}
\frac{\di\rho_i}{\di t}(T)>0.
\end{align*}
Because the continuity of $\frac{\di\rho_i}{\di t}$ w.r.t time $t$, we obtain that there exists a time $T_1<T$ such that
\begin{align*}
\frac{\di\rho_i}{\di t}(t)>0,\quad\forall t\in[T_1,T].
\end{align*}
This is contradicted with the fact 
\begin{align*}
\rho_i(T_1)>\rho_i(T)=0.
\end{align*}
Hence, we obtain that the solution $\bm{\rho}(t)$ will never reach on the boundary $\partial\mathcal{P}(G)$ at any time $T<\infty$. 
\end{proof}

\begin{proposition}\label{pro:convergence}
Let \( G = (V, E, \bm{\pi}) \) be a connected, locally uniformly finite and stochastically complete graph, $r\in [2,\infty]$. Then, the unique global solution $\bm{\rho}$ converges to the Gibbs density \(\bm{\rho}^* = (\rho_i^*)_{i=1}^\infty\in\mathcal{P}_0^*(G)\) under the $\ell^{r}(V,\bm{\pi})$ norm.
\end{proposition}
\begin{proof}
For any fixed time $t\in (0,T]$ and any $0<T<\infty$. Let $N_1$ be the subset of $\mathbb{N}$ such that
\begin{align*}
\rho_i(t)-\rho_i^*>0,\quad\textrm{if and only if}\quad i\in N_1.
\end{align*}
Let $N_2$ be the subset of $\mathbb{N}$ such that
\begin{align*}
\rho_i(t)-\rho_i^*<0,\quad\textrm{if and only if}\quad i\in N_2.
\end{align*}
Let $N_3$ be the subset of $\mathbb{N}$ such that
\begin{align*}
\rho_i(t)-\rho_i^*=0,\quad\textrm{if and only if}\quad i\in N_3.
\end{align*}
Notice that $\Psi_i=-\log{\left(\frac{1}{K}e^{-\Psi_i}\right)}-\log K=-\log\rho_i^*-\log K$.
Then, by the expression of Fokker-Planck equation \eqref{eq:FPE}, we have 
\begin{equation}\label{eq:FPETrans}
\begin{aligned}
\frac{\di}{\di t}(\rho_i-\rho_i^*)=\sum_{j\in N(i)}\frac{\omega_{ij}}{\pi_i}\left[(\log\rho_j-\log\rho_j^*)-(\log\rho_i-\log\rho_i^*)\right]\hat{\rho}(x_i,x_j).
\end{aligned}
\end{equation}
Moreover, for any $r\in\mathbb{N}$, because $\frac{1}{2r}\pi_i\frac{\di}{\di t}(\rho_i-\rho_i^*)^{2r}=\pi_i(\rho_i-\rho_i^*)^{2r-1}\frac{\di}{\di t}(\rho_i-\rho_i^*)$, the symmetry and absolute convergence of the above equation \eqref{eq:FPETrans}, we have
\begin{equation*}
\begin{aligned}
\frac{1}{2r}\sum_{i\in N_1}\pi_i\frac{\di}{\di t}(\rho_i-\rho_i^*)^{2r}\big|_{t}=I_1+I_2+I_3,
\end{aligned}
\end{equation*}
with
\begin{equation*}
\begin{aligned}
I_1=&\sum_{i\sim j\atop i,j\in N_1}\left\{\omega_{ij}\left[(\log\rho_j(t)-\log\rho_j^*)-(\log\rho_i(t)-\log\rho_i^*)\right](\rho_i(t)-\rho_i^*)^{2r-1}\hat{\rho}(x_i,x_j)(t)\right.\\
&\qquad\qquad\left.+\omega_{ij}\left[(\log\rho_i(t)-\log\rho_i^*)-(\log\rho_j(t)-\log\rho_j^*)\right](\rho_j(t)-\rho_j^*)^{2r-1}\hat{\rho}(x_i,x_j))(t)\right\}\\
=&\sum_{i\sim j\atop i,j\in N_1}\omega_{ij}\left[(\log\rho_j(t)-\log\rho_j^*)-(\log\rho_i(t)-\log\rho_i^*)\right]\\
&\qquad\qquad\qquad\qquad\qquad\qquad\qquad\quad\times\left[(\rho_i(t)-\rho_i^*)^{2r-1}-(\rho_j(t)-\rho_j^*)^{2r-1}\right]\hat{\rho}(x_i,x_j)(t),
\end{aligned}
\end{equation*}
and
\begin{equation*}
\begin{aligned}
I_2=\sum_{i\in N_1}\sum_{j\in N(i)\cap N_2}\omega_{ij}\left[(\log\rho_j(t)-\log\rho_j^*)-(\log\rho_i(t)-\log\rho_i^*)\right](\rho_i(t)-\rho_i^*)^{2r-1}\hat{\rho}(x_i,x_j)(t).
\end{aligned}
\end{equation*}
For $I_3$, we have
\begin{equation*}
\begin{aligned}
I_3=&\sum_{i\in N_1}\sum_{j\in N(i)\cap N_3}\omega_{ij}\left[(\log\rho_j(t)-\log\rho_j^*)-(\log\rho_i(t)-\log\rho_i^*)\right](\rho_i(t)-\rho_i^*)^{2r-1}\hat{\rho}(x_i,x_j)(t)\\
=&-\sum_{i\in N_1}\sum_{j\in N(i)\cap N_3}\omega_{ij}\left(\log\rho_i(t)-\log\rho_i^*\right)(\rho_i(t)-\rho_i^*)^{2r-1}\hat{\rho}(x_i,x_j)(t),
\end{aligned}
\end{equation*}
where we used the definition of $N_3$. It is not hard to prove $I_1<0$ if and only if $|N_1|\geq 2$, and 
$I_k<0$ if and only if $N_1\neq\emptyset$ and $N_k\neq\emptyset$ for $k=2,3$. Hence, we have 
\begin{align}\label{inq:important1}
\frac{1}{2r}\sum_{i\in N_1}\pi_i\frac{\di}{\di t}(\rho_i-\rho_i^*)^{2r}<0,\quad\textrm{if}~N_1\neq\emptyset.
\end{align}
Similar calculation implies
\begin{align}\label{inq:important2}
\frac{1}{2r}\sum_{i\in N_2}\pi_i\frac{\di}{\di t}(\rho_i-\rho_i^*)^{2r}<0,\quad\textrm{if}~N_2\neq\emptyset,
\end{align}
and
\begin{align}\label{inq:important3}
\frac{1}{2r}\sum_{i\in N_3}\pi_i\frac{\di}{\di t}(\rho_i-\rho_i^*)^{2r}=0, \quad\textrm{if}~N_3\neq\emptyset.
\end{align}
Notice that if $N_1\cup N_2=\emptyset$, there is nothing to prove. Next, we assume that $N_1\cup N_2\neq\emptyset$. Combining inequalities \eqref{inq:important1}, \eqref{inq:important2} and \eqref{inq:important3}, we have
\begin{equation}\label{inq:decay}
\begin{aligned}
&\frac{1}{2r}\frac{\di}{\di t}\left\|\bm{\rho}-\bm{\rho}^*\right\|_{\ell^{2r}(V,\bm{\pi})}^{2r}\\
=&\frac{1}{2r}\sum_{i=1}^\infty\pi_i\frac{\di}{\di t}(\rho_i-\rho_i^*)^{2r}\\
=&\frac{1}{2r}\sum_{i\in N_1}\pi_i\frac{\di}{\di t}(\rho_i-\rho_i^*)^{2r}+\frac{1}{2r}\sum_{i\in N_2}\pi_i\frac{\di}{\di t}(\rho_i-\rho_i^*)^{2r}+\frac{1}{2r}\sum_{i\in N_3}\pi_i\frac{\di}{\di t}(\rho_i-\rho_i^*)^{2r}<0,
\end{aligned}
\end{equation}
for any $\bm{\rho}\neq\bm{\rho}^*$, $t\in [0,T]$ and $0<T<\infty$. Since $\bm{\rho}-\bm{\rho}^*\in\ell^{2r}(V,\bm{\pi})$, $\frac{\di}{\di t}(\bm{\rho}-\bm{\rho}^*)\in\ell^{2r}(V,\bm{\pi})$ and $\bm{\rho}-\bm{\rho}^*$ is continuously differentiable in $\ell^{2r}(V,\bm{\pi})$, take Fr\'{e}chet derivative to $\left\|\bm{\rho}-\bm{\rho}^*\right\|_{\ell^{2r}(V,\bm{\pi})}^{2r}$ and use chain rule, we can obtain the first equality holds in \eqref{inq:decay}. Moreover, $\frac{1}{2r}\frac{\di}{\di t}\left\|\bm{\rho}-\bm{\rho}^*\right\|_{\ell^{2r}(V,\bm{\pi})}^{2r}<0$ holds for $t\in [0,T]$ and any $T>0$ shows it holds for $t\in [0,\infty)$. We also obtain that $\frac{1}{2r}\frac{\di}{\di t}\left\|\bm{\rho}-\bm{\rho}^*\right\|_{\ell^{2r}(V,\bm{\pi})}^{2r}<0$ if $\left\|\bm{\rho}-\bm{\rho}^*\right\|_{\ell^{2r}(V,\bm{\pi})}>0$, and $\frac{1}{2r}\frac{\di}{\di t}\left\|\bm{\rho}-\bm{\rho}^*\right\|_{\ell^{2r}(V,\bm{\pi})}^{2r}\to 0$ if and only if $\left\|\bm{\rho}-\bm{\rho}^*\right\|_{\ell^{2r}(V,\bm{\pi})}\to0$. Hence, the solution $\bm{\rho}$ convergent to $\bm{\rho}^*$ under the $\ell^{2r}(V,\bm{\pi})$ norm for $r\in\mathbb{N}$, by an easily contradiction argument.

For the case $r\in[2,\infty]/2\mathbb{Z_+}$, we use Sobolev embedding theorem to obtain the corresponding results, where $2\mathbb{Z}_+$ is the positive even number set.
\end{proof}

Combining the results of Proposition \ref{pro:gradientflow}, \ref{pro:existence} and \ref{pro:convergence}, we finished the proof of Theorem \ref{thm:globalexistence}.

\quad\\
\noindent\textbf{Data availability.} No data was used for the research described in the article.\\
\\
\noindent\textbf{Conflict of interest.} The authors declare that they have no conflict of interest.\\
\\
\noindent\textbf{Acknowledgements.} There is no acknowledgement.

\appendix

\bibliographystyle{plain}
\bibliography{OP_FP_G}

\end{document}